\documentclass{article}
\usepackage{amsmath,amsthm,amsfonts}
\newtheorem{theorem}{Theorem}
\newtheorem{lemma}{Lemma}
\newtheorem {corollary}{Corollary}
\newtheorem {proposition}{Proposition}

\theoremstyle{definition}
\newtheorem{definition}{Definition}
\newtheorem {remark}{Remark}

\DeclareSymbolFont{rsfscript}{OMS}{rsfs}{m}{n}
\DeclareSymbolFontAlphabet{\mathrsfs}{rsfscript}

\DeclareMathAlphabet{\mathbbold}{U}{bbold}{m}{n}

\def\k{\mathbbold{k}}

\DeclareFontFamily{OMS}{rsfs}{\skewchar\font'177}
\DeclareFontShape{OMS}{rsfs}{m}{n}{%
      <5> rsfs5
      <6> <7> rsfs7
      <8> <9> <10> rsfs10
      <10.95> <12> <14.4> <17.28> <20.74> <24.88> rsfs10
      }{}
\def\calA{\mathrsfs{A}}
\def\calB{\mathrsfs{B}}
\def\calC{\mathrsfs{C}}
\def\calD{\mathrsfs{D}}
\def\calE{\mathrsfs{E}}
\def\calF{\mathrsfs{F}}
\def\calG{\mathrsfs{G}}

\def\calO{\mathrsfs{O}}

\def\calR{\mathrsfs{R}}
\def\calS{\mathrsfs{S}}

\def\calV{\mathrsfs{V}}

\newcommand{\Com}{{\mathrsfs{C}\!om}}
\newcommand{\PP}{{\mathrsfs{P}_2}}
\newcommand{\LIE}{\mathrsfs{L}\!ie_2}

\newcommand{\ASS}{{\ensuremath{\mathrsfs{A}\!ss_2}}}
\newcommand{\LL}{{\ensuremath{\mathrsfs{L}\!\mathrsfs{L}_{\hbar}}}}
\newcommand{\LLL}{{\ensuremath{\mathrsfs{L}\!\mathrsfs{L}_{2,\hbar_1,\hbar_2}}}}
\newcommand{\GV}{{\ensuremath{\mathrsfs{G}\!\mathrsfs{V}_{2,\hbar_1,\hbar_2}}}}

\DeclareMathOperator{\inv}{inv}

\makeatletter
\let\@newpf\proof \let\proof\relax 
\newenvironment{proof}{\@newpf[\proofname]}{\qed\endtrivlist}
\def\hm#1{#1\nobreak\discretionary{}{\hbox{\m@th$#1$}}{}}
\makeatother

\begin{document}

%  Headings
%
%\renewcommand{\evenhead}{Vladimir Dotsenko}
%\renewcommand{\oddhead}{An operadic approach to deformation quantization of compatible Poisson brackets, I}

%  Titlepage
%
\thispagestyle{empty}

%\FirstPageHead{*}{*}{20**}{\pageref{firstpage}--\pageref{lastpage}}
%  Parameters: Volume, number, year, page range, paper type
%  'Article' could be changed to 'Letter' or 'Review Article'

\title{An operadic approach\\ to deformation quantization\\ of compatible Poisson brackets}

%\label{firstpage}

\author{Vladimir Dotsenko}

\maketitle 
\bibliographystyle{bibstand}

\begin{abstract}
An analogue of the Livernet--Loday operad for two compatible brackets, which is a flat deformation of the bi-Hamiltonian operad is constructed. The Livernet--Loday operad can be used to define $\star$-products and deformation quantization for Poisson structures. The constructed operad is used in the same way, introducing a definition of operadic deformation quantization of compatible Poisson structures.
\end{abstract}

\section{Introduction}

Throughout the text, algebras and operads are defined over an arbitrary field $\k$ of characteristic zero (unless the ground field or ring is specified explicitly).

The general deformation quantization problem for Poisson structures is set, for example, by Bayen--Flato--Fronsdal--Lichnerowicz--Stern\-heimer \cite{Lichn}. Let us give a pure algebraic definition. 

\begin{definition}A $\star$-product on a vector space $A$ is a $\k[[\hbar]]$-linear associative product on $A[[\hbar]]$ such that the  algebra $A=A[[\hbar]]/\hbar A[[\hbar]]$ is commutative. Thus, for $a,b\in A$
\[
a\star b=a*_0b+\hbar a*_1b+\ldots+\hbar^ka*_kb+\ldots\quad \text{with}\quad a*_kb\in A, \quad k\ge0
\]
\end{definition}
It is immediate to check that $A$ with the operations $a\cdot b=a*_0b$ and $\{a,b\}=a*_1b-b*_1a$ becomes a Poisson algebra, i.e. the product $a\cdot b$ is commutative and associative, and the bracket $\{a,b\}$ is a Lie bracket satisfying the Leibniz rule. 
\begin{definition}
For a Poisson algebra $A$, its deformation quantization is a $\star$-product on $A$ such that the associated Poisson algebra structure coincides with the original structure on $A$. 
\end{definition}

Some constructions of deformation quantization are known now for the case which was the most important for \cite{Lichn}, namely the algebra of functions on a smooth Poisson manifold; see, for example, the work of Kontsevich \cite{Ko}. It seems to be much more difficult to deal with deformation quantization of two compatible Poisson brackets~-- even on the level of introducing the problem and giving the necessary definitions. While for a Poisson bracket the quantum object should possess a simple natural structure (namely, an associative algebra structure), for two compatible brackets it is not clear at a glance which structures in principle one should expect. For example, an obvious idea that a pencil of Poisson structures can be quantized into a pencil of associative products is wrong. Note that although one can apply the Kontsevich's deformation quantization techniques to get a 2-parametric family of associative products, these products do not form a pencil, i.e. are not obtained by taking linear combinations of two associative products, and in general the algebraic relations between the different products from this family are unknown.

Using our previous results on the bi-Hamiltonian operad from \cite{BDK,DK}, we generalise the Livernet--Loday idea of an operadic definition of deformation quantization from the Poisson structures to the bi--Hamiltonian structures. According to our computation, the following algebraic structure should be related to the corresponding deformation quantization problem.

\begin{definition}
\label{maindef}
The operad $\ASS$ is a quadratic operad generated by two binary operations $\star\colon a,b\mapsto a\star b$ and
$[\cdot,\cdot]:a,b\mapsto[a,b]$ such that the following identities hold for each elements of any algebra over this operad:
\begin{itemize}
\itemsep-2pt
 \item $a\star (b\star c)=(a\star b)\star c$,
 \item $[a,b]+[b,a]=0$,
 \item $[a,[b,c]]+[b,[c,a]]+[c,[a,b]]=0$,
 \item $[a,b\star c+c\star b]=[a,b]\star c+b\star [a,c]+[a,c]\star b+c\star[a,b]$,
 \item $[a,b\star c-c\star b]+[b,c\star a-a\star c]+[c,a\star b-b\star a]+$\\
  \phantom{a}\hfill $+a\star [b,c]-[b,c]\star a+b\star [c,a]-[c,a]\star b+c\star [a,b]-[a,b]\star c=0$.
\end{itemize}
Thus the first operation is an associative product, and the second one is a Lie bracket which is compatible with the skew-symmetrization of the product and satisfies the Leibniz rule with the symmetrization of the product.
\end{definition}

The paper is organized as follows. We begin with some known facts on the operad defined by Livernet and Loday which can be used to define $\star$-products operadically. We prove that this operad is related to the algebras of locally constant functions on the complements of hyperplane arrangements \cite{GV}. Then we define a family of associative algebras (depending on a point of~$\k^2$) which provides a deformation for dual spaces of components of the bi-Hamiltonian operad (which are double even Orlik--Solomon algebras from \cite{BDK}); we call them double Gelfand--Varchenko algebras. The collection of these algebras possesses a cooperad structure, and it turns out that the dual operad is isomorphic to $\ASS$ everywhere except for the origin. We prove that this operad is Koszul and that our deformation of the bi-Hamiltonian operad is flat. (The proof uses the theorem on filtered operads proved in \cite{Khor}.) Finally, we discuss some aspects of the operadic problem of deformation quantization for compatible Poisson brackets. We hope to discuss examples of such a quantization and further results on it elsewhere.

Throughout the paper we use quadratic operads and the Koszul duality for operads. The corresponding definitions and notation can be found in \cite{DK,GK,MSS}.

\section{The Livernet--Loday operad and related algebras}

\begin{definition}
The double Livernet--Loday operad $\LL$ is a quadratic operad over $\k[[\hbar]]$ generated by a commutative operation $a,b\mapsto a\cdot b$ and a skew-commutative operation $a,b\mapsto[a,b]$ satisfying the identities
\begin{itemize}
\itemsep-2pt
 \item[(i)] $[a\cdot b,c]=a\cdot[b,c]+[a,c]\cdot b,$
 \item[(ii)] $[a,[b,c]]+[b,[c,a]]+[c,[a,b]]=0,$ 
 \item[(iii)] $(a\cdot b)\cdot c-a\cdot(b\cdot c)=\hbar^2[b,[a,c]]$.
\end{itemize}
\end{definition}

\begin{proposition}[\cite{MR}]
\label{ll}
For each $\hbar\ne0$ the operad $\LL$ is isomorphic to the associative operad. If $\hbar=0$, then the operad $\LL$ is isomorphic to the Poisson operad.
\end{proposition}

The relevance of this operad for deformation quantization is explained by the following proposition.

\begin{proposition}[\cite{MR}]\label{starprod}
A $\star$-product on a vector space $V$ is the same as an $\LL$-algebra structure on the $\k[[\hbar]]$-module $V[[\hbar]]$. 
\end{proposition}

The following is clear from the proof of the previous statement in \cite{MR}.

\begin{proposition}\label{defq}
A deformation quantization of a Poisson algebra $A$ is a $\star$-product on this algebra such that the operations on the vector space $A=A[[\hbar]]/\hbar A[[\hbar]]$ induced from the $\LL$-structure coincide with the original Poisson operations on~$A$.
\end{proposition}

Markl and Remm \cite{MR} proved that the operad $\LL$ is a Hopf operad. It follows that the dual spaces to its components
are associative algebras. The following fact is quite surprising though.

\begin{theorem}
The associative algebra dual to the $n$th component of the Livernet--Loday operad can be presented by  
generators $x_{ij}$, $1\le i,j\le n$, and relations
\begin{gather*}
x_{ij}+x_{ji}=\hbar\\
x_{ij}x_{jk}+\hbar x_{ik}=x_{ij}x_{ik}+x_{jk}x_{ik}\\
x_{ij}^2=\hbar x_{ij}
\end{gather*}
Thus for $\hbar\ne0$ it is isomorphic to the Gelfand--Varchenko algebra of locally constant functions on the complement of the real hyperplane arrangement $A_{n-1}$ (studied in \cite{GV}).
\end{theorem}

\begin{proof}
A straightforward calculation shows that the above relations are satisfied in the duals of the components. Thus it remains to prove that the mappings from the corresponding abstract algebras to the duals $\LL(n)^*$ are isomorphisms. It is easy to see that these abstract algebras are isomorphic to the Gelfand--Varchenko algebras. Thus the dimensions of these algebras are indeed equal to the dimensions of the duals (the number of regions for the arrangement $A_{n-1}$ is equal to $n!$, which is also the dimension of the $n$th component of the associative operad), and it is enough to check the surjectivity, which is clear.
\end{proof}

\section{Double Gelfand--Varchenko algebras}

Although the operad $\ASS$ can be (and, in fact, is) studied separately from the double Gelfand--Varchenko algebras, the construction of these algebras was crucial in our work and so we discuss them here. Moreover, one can easily check that~--- under reasonable restrictions~--- these give the essentially unique deformation of double even Orlik--Solomon algebras and thus lead to a distinguished approach to deformation quantization.

\begin{definition}
The double Gelfand--Varchenko algebra $\GV(n)$ is an associative commutative algebra with generators $x_{ij}$ and $y_{ij}$, $1\le i,j\le n$, and relations
\begin{gather*}
x_{ij}+x_{ji}=\hbar_1, y_{ij}+y_{ji}=\hbar_2\\
x_{ij}x_{jk}+\hbar_1x_{ik}=x_{ij}x_{ik}+x_{jk}x_{ik}\\
x_{ij}y_{jk}+y_{ij}x_{jk}+\hbar_2x_{ik}+\hbar_1y_{ik}=x_{ij}y_{ik}
+y_{ij}x_{ik}+x_{jk}y_{ik}+y_{jk}x_{ik}\\
y_{ij}y_{jk}+\hbar_2y_{ik}=y_{ij}y_{ik}+y_{jk}y_{ik}\\
x_{ij}^2=\hbar_1x_{ij}, y_{ij}^2=\hbar_2y_{ij}, 2x_{ij}y_{ij}=\hbar_2x_{ij}+\hbar_1y_{ij}
\end{gather*}
\end{definition}

It turns out that the collection of double Gelfand--Varchenko algebras possesses a cooperad structure. The corresponding formulae coincide with the formulae from \cite{BDK}.\footnote{As usual, we define the algebra $\GV(I)$ for an arbitrary finite set~$I$ (indexing the generators); thus, $\GV(n)=\GV(\{1,\ldots,n\})$, and cocompositions are some mappings from $\GV(I\sqcup J)$ to $\GV(I\sqcup\{*\})\otimes \GV(J)$.} (One of the restrictions fixing the deformation which were mentioned above, is that the cocomposition is given by the same formulae. It looks too restrictive at a glance, but it is true in the case of the Poisson operad.) Namely, we have the following
\begin{definition}
Let $I$, $J$ be finite sets and $*\notin I$. Let us define an algebra homomorphism 
\[
\rho_{IJ}\colon \GV(I\sqcup J)\to \GV(I\sqcup\{*\})\otimes \GV(J)
\]
by
\[
\rho_{IJ} (\square_{ij})=
\left\{
\begin{aligned}
&\square_{ij}\otimes 1, {\rm if\ }i,j\in I,\\
&1\otimes \square_{ij}, {\rm if\ }i,j\in J,\\
&\square_{i*}\otimes 1, {\rm if\ }i\in I, j\in J,\\
&\square_{*i}\otimes 1, {\rm if\ }i\in J, j\in I,
\end{aligned}
\right.
\]
where $\square$ stands for either of the letters $x$, $y$.
\end{definition}

\begin{lemma}
For all $I,J$ the homomorphism $\rho_{IJ}$ is well-defined.
\end{lemma}

It is easy to check that these mappings satisfy all the cocomposition axioms, and so define a cooperadic structure. To make the computations in the dual operad more simple and transparent, we used the following technical result.

\begin{proposition}
The algebra $\GV(n)$ is isomorphic to the associative commutative algebra with generators $a_{ij}$, $b_{ij}$, $1\le i,j\le n$ and the relations
\begin{gather*}
a_{ij}+a_{ji}=b_{ij}+b_{ji}=0\\
a_{ij}a_{jk}+a_{jk}a_{ki}+a_{ki}a_{ij}=\hbar_1^2\\
a_{ij}b_{jk}+b_{ij}a_{jk}+a_{jk}b_{ki}+b_{jk}a_{ki}+a_{ki}b_{ij}+b_{ki}a_{ij}=2\hbar_1\hbar_2\\
b_{ij}b_{jk}+b_{jk}b_{ki}+b_{ki}b_{ij}=\hbar_2^2\\
a_{ij}^2=\hbar_1^2, b_{ij}^2=\hbar_2^2, a_{ij}b_{ij}=\hbar_1\hbar_2
\end{gather*}
\end{proposition}

\section{The double Livernet--Loday operad}

\begin{definition}
The double Livernet--Loday operad $\LLL$ is generated over $\k[[\hbar_1,\hbar_2]]$ by a commutative operation $a,b\mapsto a\cdot b$ and two skew-commutative operations $a,b\mapsto \{a,b\}$ and $a,b\mapsto[a,b]$ satisfying the identities
\begin{itemize}
\itemsep-2pt
\item[(i)] $[a\cdot b,c]=a\cdot[b,c]+[a,c]\cdot b,$
\item[(ii)] $\{a\cdot b,c\}=a\cdot\{b,c\}+\{a,c\}\cdot b,$
\item[(iii)] $\{a,\{b,c\}\}+\{b,\{c,a\}\}+\{c,\{a,b\}\}=0,$
\item[(iv)] $[a,[b,c]]+[b,[c,a]]+[c,[a,b]]=0,$ 
\item[(v)] $\{a,[b,c]\}+\{b,[c,a]\}+\{c,[a,b]\}+[a,\{b,c\}]+[b,\{c,a\}]+[c,\{a,b\}]=0,$
\item[(vi)] $(a\cdot b)\cdot c-a\cdot(b\cdot c)=\hbar_1^2[b,[a,c]] +\hbar_1\hbar_2([b,\{a,c\}] +\{b,[a,c]\})+\hbar_2^2\{b,\{a,c\}\}$.
\end{itemize}
\end{definition}

\begin{proposition}
For each point $(\hbar_1,\hbar_2)$ except for the origin the operad $\LLL$ is isomorphic to the operad $\ASS$. If $(\hbar_1,\hbar_2)=(0,0)$ then the operad $\LLL$ is isomorphic to the bi-Hamiltonian operad.
\end{proposition}

\begin{proof}
The case of the origin is clear. Let $(\hbar_1,\hbar_2)\ne(0,0)$. Consider a new operation 
\[
a,b\mapsto \hbar_1[a,b]+\hbar_2\{a,b\}
\]
It is clear that the operad generated by this operation and $\cdot$ is isomorphic
to $\LL$ (and thus, according to Proposition \ref{ll}, to the associative operad). The remaining bracket gives the second operation of $\ASS$.
\end{proof}

For the sake of completeness, we describe here the quadratic dual operad. A straightforward calculation proves the following

\begin{proposition}
The quadratic dual operad $\LLL^!$ is generated over $\k[[\hbar_1,\hbar_2]]$ by a skew-commutative operation $a,b\mapsto \{a{,}b\}$ and two commutative operations $a,b\mapsto a\star_1 b$ and $a,b\mapsto a\star_2 b$ satisfying the identities
\begin{itemize}
\itemsep-2pt
\item[(i)] $\{a\star_1 b,c\}=a\star_1\{b,c\}+\{a,c\}\star_1 b$,
\item[(ii)] $\{a\star_2 b,c\}=a\star_2\{b,c\}+\{a,c\}\star_2 b$,
\item[(iii)] $\{a,\{b,c\}\}+\{b,\{c,a\}\}+\{c,\{a,b\}\}=0,$
\item[(iv)] $(a\star_i b)\star_j c-a\star_i(b\star_j c)=\hbar_i\hbar_j\{b,\{a,c\}\}\ (i,j\in\{1,2\})$, 
\item[(v)] $(a\star_i b)\star_j c=(a\star_j b)\star_i c\quad (i,j\in\{1,2\})$. 
\end{itemize}
\end{proposition}

\section{Filtrations and Koszulness}
\label{sec_filtration}

Here we recall the results of A.~Khoroshkin that are crucial for proving the Koszulness of operads in several important cases.

Let $\Gamma$ be a set-theoretic operad. It is called an ordered operad if each set $\Gamma(n)$ is ordered and this order is compatible with operadic compositions: for $a<b\in\Gamma(n)$, $c\in\Gamma(k)$, $i=1,\ldots,k$, $j=1,\ldots,n$, we have
\[
c\circ_i a< c\circ_i b, \quad a\circ_j c< b\circ_j c
\]
Suppose that for each $n$ and each $\alpha\in\Gamma(n)$ we have a subspace $F_{\alpha}\subset \calO(n)$. The collection of these subspace is called a $\Gamma$-valued filtration on $\calO$, if
\begin{itemize}
\itemsep-2pt
\item 
for all $\alpha\leq \beta$ we have $F_{\alpha}\subset F_{\beta}$,
\item 
for the maximal element $\delta_n\in\Gamma(n)$ we have $F_{\delta_n}=\calO(n)$, 
\item 
fompositions are compatible with the filtration
\[
\gamma_{l,m_1,...,m_l}\left(F_{\alpha_l}\otimes F_{\beta_1}\otimes\ldots\otimes F_{\beta_l}\right)
\subset F_{\alpha_l(\beta_1,...,\beta_l)}
\]
\end{itemize}

The associated graded operad $gr^{F}\calO(n)$ is, by the definition, the sum
\[
\bigoplus_{\alpha\in\Gamma(n)} \frac{F_{\alpha}}{\Sigma_{\alpha'<\alpha} F_{\alpha'}}
\]

A $\Gamma$-valued filtration on $\calO$ is \emph{generated by binary operations} if for each $\alpha$ the space $F_{\alpha}$ coincides with the span of compositions of all binary operations that belong to some $F_{\alpha'}$ with $\alpha'\le\alpha$.

The following theorem is proved in \cite{Khor}; both theorems with proofs generalise Theorem 7.1 of \cite{PP} from the quadratic algebras to quadratic operads.

\begin{theorem}
\label{filtrations}
Suppose that a $\Gamma$-valued filtration on a quadratic operad $\calA$ (with generators $\calO_\calA$ and relations $\calR_\calA$) is generated by binary operations and satisfies the following conditions:
\begin{itemize}
\itemsep-2pt
\item[(i)] the operad $gr^{F}\calA$ does not have nontrivial relations of degree $3$ (between the elements belonging to $gr^{F}\calA(4)$),
\item[(ii)] the quadratic operad $qgr^{F}\calA$ with generators $gr^{F}\calO_\calA$ and relations $gr^{F}\calR_\calA$ is Koszul.
\end{itemize}
Then the operad $gr^{F}\calA$ is quadratic (and thus isomorphic to $qgr^F(\calA)$), and the operad $\calA$ is Koszul.
\end{theorem}

As an example we give a proof of the theorem on distributive laws between Koszul operads \cite{M}. The proof of this statement in \cite{M} is partially based on an erroneous formula for the composition in the symmetric case \cite[Prop.~1.7]{M} and is therefore incomplete. The proof below is essentially due to A.~Khoroshkin \cite{Khor}.
  
Let $\calA$ and $\calB$ be two quadratic Koszul operads. Denote by $\calO_\calA$, $\calO_\calB$ the corresponding sets of generators and by $\calR_\calA$, $\calR_\calB$ the spaces of relations. Denote also by $\calO_\calA\bullet\calO_\calB$ the subspace in the free operad $\calF_{\calO_\calA\oplus\calO_\calB}$ generated by all elements $\phi(1,\psi)$ with $\phi\in\calO_\calA$, $\psi\in\calO_\calB$ (and $1$ stands for the identity unary operation). The notation $\calO_\calB\bullet\calO_\calA$ has the analogous meaning. Assume that there is an $S_3$-equivariant mapping 
 \[
d\colon\calO_\calB\bullet\calO_\calA\rightarrow\calO_\calA\bullet\calO_\calB
 \]
Let $\calC$ be the quadratic operad with generators $\calO_\calC=\calO_\calA\oplus\calO_\calB$ and relations
$\calR_\calC=\calR_\calA\oplus\calD\oplus\calR_\calB$, where 
 \[
\calD=\{x-d(x)\mid x\in\calO_\calB\bullet\calO_\calA\}
 \]

\begin{theorem}[\cite{M,Khor}]
Assume that the natural mapping
 \[
\xi\colon(\calA\circ\calB)(4)\to\calC(4)
 \] 
is an isomorphism. Then the operad $\calC$ is Koszul, and $\mathbb{S}$-modules $\calA\circ\calB$ and $\calC$ are isomorphic.
\end{theorem}

\begin{proof}

The main ingredient of the proof is some particular set-theoretic operad $\Gamma$ and a $\Gamma$-filtration on $\calC$ for which the associated graded operad is isomorphic to $\calA\circ\calB$ (where the composition in the reverse order is equal to zero). 

Consider the free set-theoretic operad generated by two binary operations $a$ and $b$. Let $T$ be an arbitrary element of this operad, i.e. a binary tree whose internal vertices are labeled by $a$ and $b$. By the definition, the degree $\deg(T)$ is equal to the number of $a$'s among the labels, and the number of inversions $\inv(T)$ is equal to the number of pairs $(v_1,v_2)$ of internal vertices of $T$ where $v_2$ belongs to the set of descendants of $v_1$, the label of $v_2$ is equal to $a$, and the label of $v_1$ is equal to $b$. We have
\[\calF_2(n) = \bigsqcup_{k=0}^{n-1} \calF_2(n)_k\ {\mbox{ with }} \calF_2(n)_k:=\{T| \deg(T)=k\}\]
Let $m_{n,k}$ be the maximal number of inversions for the elements of $\calF_2(n)_k$.

Let us define the operad of inversions $\calF_2^{\inv}$. Its component $\calF_2^{\inv}(n)$ is a disjoint union of sets
 \[
\bigsqcup_{k=0}^{n-1} \{s_0^{k}(n),\ldots,s_{m_{n,k}}^{k}(n)\} 
 \]
the ordering of this component is lexicographic (to compare $s_i^j(n)$ and $s_k^l(n)$, we use the lexicographic ordering for the pairs $(j,i)$ and  $(l,k)$), and the compositions are given by
  \[
s_i^{l}(n)(s_{j_1}^{k_1}(r_1),\ldots,s_{j_n}^{k_n}(r_n))
:=s_{\max(\inv(\alpha(\beta_1,\ldots,\beta_n)))}^{l+k_1+\ldots+k_n}(r_1+\ldots+r_n)
  \]
where maximum is taken over all collections $(\alpha,\beta_1,\ldots,\beta_n)$ with 
$$
\inv(\alpha)=i,\ \inv(\beta_1)=j_1,\ldots,\inv(\beta_n)=j_n,\
\deg(\alpha)=l,\ 
\deg(\beta_1)=k_1,\ldots,\deg(\beta_n)=k_n
$$

There exists a natural $\calF_2^{\inv}$-valued filtration $F$ on $\calC$: it assigns the label $a$ to the generators belonging to $\calA$,
and the label $b$ to the generators belonging to $\calB$. The associated graded operad $qgr^F\calC$ is isomorphic to $\calA\circ\calB$ (with zero reverse composition). Indeed, the associate graded relations include $\calR_\calA$ and $\calR_\calB$ remaining as is, while the relations  $x-d(x)$ are transformed into $x$, since $d(x)$ belongs to the previous level of the filtration with respect to the number of inversions.

The operad $qgr^F\calC$ is Koszul, since for the corresponding Koszul complexes we have the following isomorphisms and quasi-isomorphisms:
\begin{multline*}
K^{\bullet}(qgr^F\calC)\simeq(\calA\circ\calB)\circ((\calA\circ\calB)^!)^{\vee}
\simeq
\calA\circ \left(\calB\circ (\calB^{!})^{\vee}\right)\circ(\calA^{!})^{\vee}
\stackrel{quis}{\longrightarrow} \calA\circ(\calA^{!})^{\vee} \stackrel{quis}{\longrightarrow}\k
\end{multline*}
Moreover, the operad $gr^F\calC$ has no nontrivial relations of degree $3$, since the corresponding components of $gr^F\calC$ and $qgr^F\calC$ are isomorphic. It follows that $\calC$ is Koszul, and $qgr^F\calC\simeq gr^F\calC$. Taking the associated graded operad does not affect the $\mathbb{S}$-module structure, so $\calC\simeq\calA\circ\calB$. 
\end{proof}

This proof can be easily generalised for the following problem.

Let again $\calA$ and $\calB$ be two quadratic operads (all notation remains the same). Assume that there are $S_3$-equivariant mappings
 \[
s\colon\calR_\calA\rightarrow\calO_\calB\bullet\calO_\calA\oplus\calO_\calA\bullet\calO_\calB\oplus\calF_{\calO_\calB}(3)
 \]
and
 \[
d\colon\calO_\calB\bullet\calO_\calA\rightarrow\calO_\calA\bullet\calO_\calB
 \]
Let $\calE$ be the quadratic operad with generators $\calO_\calE=\calO_\calA\oplus\calO_\calB$ and relations
$\calR_\calE=\calS\oplus\calD\oplus\calR_\calB$, where 
 \[
\calS=\{x-s(x)\mid x\in\calR_\calA\},\quad \calD=\{x-d(x)\mid x\in\calO_\calB\bullet\calO_\calA\}
 \]

\begin{theorem}
Assume that the natural mapping
 \[
\xi\colon(\calA\circ\calB)(4)\to\calE(4)
 \] 
is an isomorphism. Then the operad $\calE$ is Koszul, and  $\mathbb{S}$-modules $\calA\circ\calB$ and $\calE$ are isomorphic.
\end{theorem}

\begin{proof}
The same filtration can be used here as well. It turns out again that $qgr^F\calE\simeq\calA\circ\calB$ (the associated graded relations included $\calR_\calB$ remaining as is, the relations $x-d(x)$ are transformed into $x$ since $d(x)$ belongs to the previous level of filtration with respect to the number of inversions, while the relations $x-s(x)$ are transformed into $x$ since $s(x)$ belongs to the previous level of filtration with respect to the degree), and the remaining proof is the same.
\end{proof}

\begin{corollary}
The operad $\LLL$ is Koszul, $\mathbb{S}$-module $\LLL$ is isomorphic to $\Com\circ\LIE$. 
\end{corollary}

\begin{proof}
We apply the theorem to $\calA=\Com$, $\calB=\LIE$. A straightforward calculation shows that there are no nontrivial relations of degree $3$ in $gr^F\LLL$. 
\end{proof}

This corollary together with the results of \cite{DK} gives the following

\begin{corollary}
$\mathbb{S}$-modules $\LLL$ and $\PP$ are isomorphic, so we have constructed a flat deformation of the bi-Hamiltonian operad.
\end{corollary}

Similarly to \cite{BDK}, we can derive some information about double Gelfand--Varchenko algebras from the results about the double Livernet--Loday operad. 

The following proposition is straightforward.
\begin{proposition}
Consider the operad $\calG\!\calV^*$ dual to the cooperad of double Gelfand--Varchenko algebras. Then
\begin{itemize}
\itemsep-2pt
 \item[(i)] the relations of $\LLL$ hold in $\calG\!\calV^*$,
 \item[(ii)] quadratic relations of $\calG\!\calV^*$ are exactly the relations of $\LLL$.
\end{itemize}
\end{proposition}

It follows that there exists a natural mapping of the operads 
 \[
\pi\colon\LLL\to\calG\!\calV^*
 \] 

\begin{theorem}
The operad $\calG\!\calV^*$ is quadratic and $\pi$ is an isomorphism. 
\end{theorem}

\begin{proof}
The proof is completely analogous to the proof in the case of double even Orlik--Solomon algebras in \cite{BDK}. 
We already proved that the operad $\LLL$ is a flat deformation of the bi-Hamiltonian operad. In fact, it is clear that the  monomials introduced in \cite{BDK} provide a basis for this operad as well. These monomials can be used to prove that the pairing between $\LLL$ and the collection of double Gelfand--Varchenko algebras is nondegenerate (where we have to check the triangularity of a certain matrix), which proves the theorem.
\end{proof}

\section{Deformation quantization: an operadic definition}

We introduce the following definitions which are motivated by Propositions \ref{starprod} and \ref{defq}.

\begin{definition}
An operadic $\star_{\hbar_1,\hbar_2}$-product on a vector space $V$ is a $\LLL$-algebra structure on the $\k[[\hbar_1,\hbar_2]]$-module $V[[\hbar_1,\hbar_2]]$.
\end{definition}

\begin{definition}
An operadic deformation quantization of a bi-Hamiltonian algebra $A$ is a $\star_{\hbar_1,\hbar_2}$-product on this algebra such that the operations on the vector space $A=A[[\hbar_1,\hbar_2]]/(\hbar_1,\hbar_2)$ induced from the $\LLL$-structure coincide with the original (bi-Hamiltonian) operations on~$A$.
\end{definition}

Similarly to Proposition \ref{starprod}, one can check that for any $\LLL$-algebra $V$ over $\k[[\hbar_1,\hbar_2]]$ the product on $V$ given by 
\[
a\star_{\hbar_1,\hbar_2}b=a\cdot b+\hbar_1\{a,b\}_1+\hbar_2\{a,b\}_2
\]
is associative. Unfortunately, this definition does not possess the most nice property of the definition for the Poisson algebras. Namely, one can prove

\begin{proposition}
There exist obstructions for the existence of an operadic deformation quantization; they are nontrivial for a generic bi-Hamiltonian algebra, and even for a generic bi-Hamiltonian structure on $\k^n$.
\end{proposition}

\begin{remark}
We expect that the deformation quantization theorem holds in a weaker form for this operadic version: for a bi-Hamiltonian algebra $A$, there should exist a structure of a strong homotopy $\LLL$-algebra \cite{GK} on $A[[\hbar_1,\hbar_2]]$.
\end{remark}

\section*{Acknowledgements}
The author is grateful to B.~Feigin, A.~Khoroshkin and L.~Rybnikov for useful discussions. 
This work is partially supported by the LIEGRITS fellowship MRTN-CT-2003-505078.

\footnotesize
\noindent
Independent University of Moscow,\\
 Bolshoj Vlasievsky per. 11, 119002 Moscow,  Russia\\
\texttt{dotsenko@mccme.ru}

\end{document}